\documentclass[12pt]{article}
\usepackage[margin=2cm]{geometry}
\usepackage{amsfonts}
\usepackage{amssymb}
\usepackage{amsthm}
\usepackage{amsmath}
\usepackage{latexsym}
\usepackage{color}

\begin{document}
\newtheorem{lemma}{Lemma}
\newtheorem{pron}{Proposition}
\newtheorem{thm}{Theorem}
\newtheorem{Corol}{Corollary}
\newtheorem{exam}{Example}
\newtheorem{defin}{Definition}
\newtheorem{remark}{Remark}
\newtheorem{property}{Property}
\newcommand{\la}{\frac{1}{\lambda}}
\newcommand{\sectemul}{\arabic{section}}
\renewcommand{\theequation}{\sectemul.\arabic{equation}}
\renewcommand{\thepron}{\sectemul.\arabic{pron}}
\renewcommand{\thelemma}{\sectemul.\arabic{lemma}}
\renewcommand{\thethm}{\sectemul.\arabic{thm}}
\renewcommand{\theCorol}{\sectemul.\arabic{Corol}}
\renewcommand{\theexam}{\sectemul.\arabic{exam}}
\renewcommand{\thedefin}{\sectemul.\arabic{defin}}
\renewcommand{\theremark}{\sectemul.\arabic{remark}}
\def\REF#1{\par\hangindent\parindent\indent\llap{#1\enspace}\ignorespaces}
\def\lo{\left}
\def\ro{\right}
\def\be{\begin{equation}}
\def\ee{\end{equation}}
\def\beq{\begin{eqnarray*}}
\def\eeq{\end{eqnarray*}}
\def\bea{\begin{eqnarray}}
\def\eea{\end{eqnarray}}
\def\d{\Delta_T}
\def\r{random walk}
\def\o{\overline}

\title{\large\bf Convolution and convolution-root properties of long-tailed distributions
}
\author{\small
Hui Xu$^{1}  $~~ Sergey Foss$^{2,3}$\thanks{
Research supported by EPSRC grant No. EP/I017054/1}~~
Yuebao Wang$^{1}$\thanks{Research supported by National Science Foundation of China, Grant
No.11071182  }
\thanks{Corresponding author.
Telephone: +86 512 67422726. Fax: +86 512 65112637. E-mail:
ybwang@suda.edu.cn}
\\
{\footnotesize\it 1. School of Mathematical Sciences, Soochow
University, Suzhou 215006, China}\\ {\footnotesize\it 2. School of MACS and Maxwell Institute, Heriot-Watt
University, Edinburgh EH14 4AS, UK}\\
{\footnotesize\it 3. Sobolev's Institute of Mathematics and Novosibirsk State University,
Novosibirsk 630090, Russia}}
\date{}

\maketitle

\begin{center}
{\noindent\small {\bf Abstract }}
\end{center}

{\small
We obtain a number of new general properties, 
related to the
closedness of the class of long-tailed
distributions under convolutions,
that are of interest themselves and may be applied in many models that
deal with ``plus'' and/or ``max'' operations on heavy-tailed random variables.
We analyse the closedness property under convolution roots for these distributions.
Namely, we introduce two classes of heavy-tailed
distributions that are not long-tailed and study their properties.
These examples help to provide further insights and, in particular, to show that the properties
to be both long-tailed and so-called ``generalised subexponential''
are not preserved under the convolution roots.
This leads to a negative answer to a conjecture of Embrechts and Goldie \cite{EmbGol1980, EmbGol1982} for the class of long-tailed and generalised subexponential distributions. In particular, our examples show that the following is possible: an infinitely divisible distribution belongs 
to both classes,
while its L$\acute{e}$vy measure is 
neither long-tailed nor generalised subexponential.

\medskip

\medskip

}

\section{Introduction}
\setcounter{thm}{0}\setcounter{Corol}{0}\setcounter{lemma}{0}\setcounter{pron}{0}
\setcounter{remark}{0}\setcounter{exam}{0}\setcounter{property}{0}

We assume $F$ to be a distribution on the real line,
with the (right) tail distribution function $\overline F (x)=1-F(x)$.
The notation $F_1*F_2$ is reserved for the convolution of two
distributions $F_1$ and $F_2$; further $F^{*n}=F* \ldots *F$ denotes the $n$-fold
convolution of $F$ with itself for $n\ge2$, and $F^{*1}=F$ and $F^{*0}$ denotes the  distribution 
degenerate
at zero.
All limits are taken as $x$ tends to infinity. For two positive
functions $f$ and $g$, the notation
$f(x)\sim g(x)$ means that $\lim f(x)/g(x)=1$; the notation $f(x)=o(g(x))$
means that $\lim f(x)/g(x)=0$; and $f(x)=O(g(x))$ means that $\limsup
f(x)/g(x)<\infty$. The indicator function $\textbf{I}(A)$ of an event $A$
takes the value $1$ if the event occurs and the value 
$0$
otherwise.

Recall that a distribution $F$ on the real line is {\it heavy-tailed}
if $\int_0^{\infty} e^{\beta y}F(dy)=\infty$ for all $\beta>0$, otherwise $F$ is
{\it light-tailed}.
A distribution $F$ is {\it long-tailed}, denoted by $F\in {\mathcal L}$, if
$\overline{F}(x+1)\sim\overline{F}(x)$. A distribution $F$ {\it
on the positive half-line} is {\it subexponential}, denoted by $F\in {\mathcal S}$, if
$\overline{F^{*2}}(x)\sim2\overline{F}(x)$.
A distribution $F$ {\it on the whole real line} is
subexponential if the distribution $F_+$ is subexponential, where $F_+(x) = F(x)\cdot
{\bf I} (x\ge 0)$ for all $x$, or, equivalently, if $F\in\mathcal{L}$ and if $\overline{F^{*2}}(x)\sim2\overline{F}(x)$. Note that both subexponentiality and long-tailedness are the
{\it tail properties}: if a distribution $F$ has such a property and
$\overline{F}(x)\sim \overline{G}(x)$, then $G$ also has this property.
It is known that any subexponential distribution is long-tailed and any long-tailed distribution is heavy-tailed.

More generally, let $\gamma\ge 0$ be fixed. A distribution $F$
on the whole real line
{\it belongs to the distribution class}
$\mathcal{L}(\gamma)$ if, for any fixed $c>0$,
$$\overline F(x-c)\sim \overline F(x)e^{\gamma c}.$$
A distribution $F$ {\it belongs to the class} $\mathcal{S}(\gamma)$
if $\int_0^{\infty}e^{\gamma y}F(dy)<\infty,\ F\in\mathcal{L}(\gamma)$ and if
$$
\overline {F^{*2}}(x)\sim 2\overline F(x)\int_{-\infty}^{\infty}e^{\gamma
y}F(dy).
$$

In particular,
${\mathcal L}= \mathcal{L}(0)$
and
$\mathcal{S} = \mathcal{S}(0)$.
Clearly, distributions from the class $\mathcal{L}(\gamma )$ are light-tailed if $\gamma >0$.
For all $\gamma\ge0$, the class
$\mathcal{L}(\gamma )\setminus \mathcal{S}(\gamma)$ is non-empty,
see, e.g.,
Pitman \cite{Pit1980}, Leslie \cite{Les1989}, Murphree \cite{Mur1989},
Kl\"{u}ppelberg and Villasenor \cite{KluVil1991}
and
Lin and Wang \cite{LinWan2012} 
for examples and further analysis. Some systems
research and application on the above mentioned distribution classes, please refer to
Embrechts et al. \cite{EmbKluMik1997}, Asmussen and Albrecher \cite{AsmAlb2010}, Foss et al. \cite{FosKorZac2013}, and so on.

Recall that the classes $\mathcal{S}$ and $\mathcal{L}$ were introduced by Chistyakov \cite{Chi1964} and, for $\gamma>0$, the class $\mathcal{S}(\gamma)$ of distributions supported by the positive half-line
was introduced and analysed by Chover et al. \cite{ChoNeyWai1973a, ChoNeyWai1973b}.
The class ${\cal L}$ is closely linked to slow variation ($F\in \mathcal{L}$ iff
$\overline{F} (\log x)$ is slowly varying).  For $\gamma>0$, the
class $\cal{L}(\gamma )$ was introduced
by Embrechts and Goldie \cite{EmbGol1980} and is linked to
regular variation.

It is known that if $F\in {\mathcal L}$ (or if $F\in {\mathcal S}$), then $F^{*n}\in {\mathcal L}$
(correspondingly $F^{*n} \in {\mathcal S}$), for any $n\ge 2$.
These results continue to hold when $F^{*n}$ is replaced by the
{\it compound distribution} 
$\sum_{n=0}^{\infty} p_n F^{*n}$ where $0\le p_n\le1$ 
for $n=0,1,\ldots$, $p_0<1$, $\sum_{n=0}^\infty p_n=1$, given than
$p_n$ decay to zero sufficiently fast as $n\to\infty$. 
In the case of subexponential distributions this is a classical result
(based
on ``Kesten's lemma''; see also \cite{DenFosKor2010} and the references therein
for modern results in this direction),
while the result for long-tailed distributions is
quite recent (\cite{Alb2008, LeiSia2012}).
Similar results hold for the class $\mathcal{S}(\gamma )$ for any $\gamma>0$.
Therefore, we may say that all these distribution classes are {\it closed under convolution}.

Embrechts et al. \cite{EmbGolVer1979} (see also \cite{EmbGol1981}) proved the converse result for subexponential distributions: if $F^{*n}\in {\mathcal S}$ for some $n\ge 2$, then $F\in {\mathcal S}$ (and, in turn, $F^{*m}\in {\mathcal S}$, for all $m
\ge 2$). They also proved an analogous result related to the 
compound distribution,
and then similar
results for the class 
$\mathcal{S}(\gamma )$ for any $\gamma >0$.    In short, one can say
that, for any $\gamma\ge0$, the class $\mathcal{S}(\gamma )$
is {\it closed under convolution roots}.

Embrechts and Goldie (see \cite{EmbGol1980}, page 245 and \cite{EmbGol1982}, page 270) formulated the conjecture that a similar converse result may hold for long-tailed distributions and, more generally,
for any class $\mathcal{L}(\gamma )$, $\gamma \ge 0$. 

{\bf Conjecture 1}
Let $\gamma\ge 0$. If there is $n\ge 2$ such that $F^{*n}\in \mathcal{L}(\gamma )$,
then also $F\in \mathcal{L}(\gamma )$.

The following two closely related
conjectures may be viewed as natural extensions of Conjecture~1
onto compound distributions and infinitely divisible distributions.

{\bf Conjecture 2.} Let $\gamma\ge0$ and $\sum_{n=0}^{\infty} p_n =1$, with $p_n\ge 0$ for all $n$ and $p_0+p_1<1$. If a compound distribution $\sum_{n=0}^\infty p_n F^{*n}$ belongs to the
class $\mathcal{L}(\gamma )$, 
then also $F\in \mathcal{L} (\gamma )$.

{\bf Conjecture 3.} Let $\gamma\ge0$. If an infinitely divisible distribution $H$ belongs to the class
$\mathcal{L}(\gamma )$, then the 
distribution generated by its L$\acute{e}$vy spectral measure belongs to the class  $\mathcal{L} (\gamma )$ too.


In this paper, we restrict our attention to the study of  the class of long-tailed
distributions, and also of its subclass consisting of the so-called
\emph{generalised subexponential} distributions.

A distribution $F$ is {\em generalised subexponential}, denoted
by $F\in{\mathcal{OS}}$, if $\overline{F}(x)>0$ for all $x$ and if
\begin{eqnarray*}
C^*(F)=\limsup\overline{F^{*2}}(x)/\overline{F}(x)<\infty.
\end{eqnarray*}
Note that (a) for any heavy-tailed distribution $F$ on the whole real line,
$C^*(F) \ge 2$ (see Theorem 1.2 in \cite{YWC2010} for this and further results);
(b) clearly, $C^*(F) \ge \liminf \overline{F^{*2}}(x)/\overline{F}(x) \ge 2$ for any
distribution on the positive half-line.

The class $\mathcal{OS}$ was first introduced
by Kl\"{u}ppelberg \cite{Klu1990} for distributions on the positive
half-line and was called ``weakly idempotent''.
Later Shimura and Watanabe \cite{ShiWat2005} called it
``O-subexponential'', or  ``generalised subexponential'', by analogy to``O-regularly varying'' in the terminology
of Bingham et al.~\cite{BGT1987}. The definition of the
class $\mathcal{OS}$ was extended in \cite{W2008} 
to the whole
real line. 

In this paper, we prove a number of novel properties of long-tailed distributions
(see Theorem \ref{th1}) that, in particular, allow us to
provide a number of counter-examples to Conjectures 1-3 (see Theorem \ref{th2} and Proposition
 \ref{classF2}) 
where the class $\mathcal{L}$ is replaced by the class $\mathcal{L}\cap\mathcal{OS}$.
We also provide a simple sufficient condition for the equivalence ``$F\in {\cal L}$ if and only if
$F^{*2}\in {\cal L}$'' to hold, see Proposition \ref{pron1}.
Similar problems for light-tailed distributions (with counterexamples to Conjectures 2-3)
will be analysed in a companion paper.

The remainder of this paper is organised as follows. In Section 2 we formulate
and discuss our
main results and their corollaries.
In 
Section 3 we prove Theorem \ref{th1} and Corollary \ref{cor001}.
The proofs of Theorem \ref{th2}, Proposition \ref{pron1} and Lemma \ref{gengen}
are given in Section 4.
Finally, the appendix includes comments related to the condition \eqref{thm101}
and  a sketch of the proof of Proposition \ref{classF2}.

\section{Main results and 
related discussions}
\setcounter{thm}{0}\setcounter{Corol}{0}\setcounter{lemma}{0}\setcounter{pron}{0}
\setcounter{remark}{0}\setcounter{exam}{0}\setcounter{property}{0}

To formulate our first result, we need further notation.
For a distribution $F$ 
and any constants $a\le b$, we let $F(a,b] = F(b)-F(a) = 
\overline{F}(a)-\overline{F}(b)$.
Let $X_1,X_2,\ldots$ be independent (not necessarily identically distributed) random variables with corresponding distributions
$F_1,F_2,\ldots$. For $n=0,1,\ldots$, let $S_n=\sum_{i=1}^n X_i$ be the partial sum
with distribution $H_n=F_1*\cdots *F_n$, where $H_0$ 
degenerates at $0$. Let $\tau$ be an independent
counting random variable with distribution function $G(x) = \sum_{n\le x} p_n$ where
$p_n={\mathbf P} (\tau =n)$, $n=0,1,\ldots$. We denote by $H_{\tau}$ the
distribution 
of the random sum $S_{\tau}=\sum_{i=1}^{\tau} X_i$. Clearly,
$H_{\tau}=\sum_{n=0}^{\infty} p_nH_n$.
In the particular case where $\{X_i,i\ge1\}$ are i.i.d. with common distribution $F$,
we have $H_n=F^{*n}$ for $n=0,1,\ldots$ and we also use notation $F^{*\tau}$
for
$H_{\tau}= \sum_{n=0}^{\infty} p_n F^{*n}$.

\begin{thm}\label{th1} (1) Let $n\ge 2$. \\
(1a) If $H_n\in\mathcal{L}$, then
\begin{equation}\label{1011}
F_i(x-c,x+c] = o(\overline{H_n}(x))
\quad \mbox{and}
\quad
H_i(x-c,x+c] = o(\overline{H_n}(x)),
\end{equation}
 for any $c>0$ and
all $i=1,\ldots, n$.\\
(1b) Assume $H_m\in\mathcal{L}$ for some $1\le m\le n$ and
\begin{equation}\label{101}
F_i(x-c,x+c]
= o(\overline{H_m}(x)),
\end{equation}
for some $c>0$ and all $i=m+1,\ldots, n$. Then
$H_n\in\mathcal{L}$. \\
(2) Let $\tau$ be an independent counting random variable with bounded support:
$\sum_{k=0}^n p_k=1$ and $p_n>0$, for some $n\ge 1$.
Then $H_{\tau}\in\mathcal{L}$ if and only if $H_n\in \mathcal{L}$.\\
(3) Assume that ${\mathbf P}(\tau \ge n)>0$ for some $n\ge 1$ and
$H_{k}\in\mathcal{L}$ for all $k\ge n$. Assume further that there exists a positive constant
$C$ such that,
for every $n=1,2,\ldots$,
the following concentration inequality holds:
\begin{equation}\label{concentration}
\sup_x H_n(x-1,x]\le C/\sqrt{n}
\end{equation}
and that, for any $\varepsilon >0$, there exists $x_0>1$ such that, for all
$k\ge n$,
\begin{equation}\label{tail}
\sup_{x\ge k(x_0-1)+x_0}
\overline{H}_k(x-1)/\overline{H}_k(x)\le 1+ \varepsilon.
\end{equation}
If, in addition, for any $a>0$,
\begin{equation}\label{102}
\overline{G}(ax) = o\left(x^{1/2}\overline{H_{n}}(x)\right),
\end{equation}
then $H_{\tau}\in \mathcal{L}$.\\
(4) Let
$F_1,F_2$ and $L_2$ be three distributions
such that $\overline{F}_2(x)\sim
\overline{L_2}(x)$ and $F_1*F_2\in\mathcal{L}$. Then
$\overline{F_1*L_2}(x)\sim
\overline{F_1*F_2}(x)$ and, therefore, $F_1*L_2\in\mathcal{L}$.
\end{thm}

\begin{remark}\label{eq-eq}
Statement (1b) of Theorem \ref{th1} is equivalent to the following:\\
Assume $F_n\in\mathcal{L}$ for some $n\ge2$ and
$$
F_i(x-c,x+c] = o\big(\overline{F_n}(x)\big),
$$
for some $c>0$ and all $i=1,\ldots, n$. Then
$H_n\in\mathcal{L}$.
\end{remark}

\begin{remark}\label{rem2000}
Condition \eqref{concentration} is very general. It holds if random variables
$X_i,\ i\ge1,$ are i.i.d. with any non-degenerate distribution (see, e.g.,
\cite{Pet1995}, Theorem 2.22).  More generally, \eqref{concentration} holds
if random variables $X_i$ are assumed to be independent, but not necessarily
identically distributed, and
there exists $c>0$ such that
\begin{equation}\label{suff}
\inf_{i\ge1} {\mathbf P} (X_i\in [-c,c])>0 \quad
\mbox{and}
\quad
\inf_{i\ge1} {\mathbf Var} (X_i \ | \ X_i\in [-c,c]) >0,
\end{equation}
see e.g. \cite{FosKor2000}, Lemma 4.1. Moreover, it is enough to assume
that \eqref{suff}
holds only for a positive proportion of the summands: if $c_n$ is the number of
$X_i$, $i\le n$ that satisfy \eqref{suff}, then $c_n/n \ge c>0$ for some $c>0$
and for all sufficiently large $n$.

Some other conditions for the concentration inequality can be found in theorems
that precede Theorem 2.22 of the book \cite{Pet1995} (e.g., Theorems 2.17 and 2.18).
\end{remark}

In the case of i.i.d. summands, Theorem \ref{th1} leads to the following corollary.

\begin{Corol}\label{cor001}
(1) Assume a distribution $F$ to be such that $F^{*n}\in \mathcal{L}$, for some $n\ge 1$.
Then $F^{*k}\in \mathcal{L}$, for all $k\ge n$.\\
(2) Let $\tau$ be a counting random variable with bounded support:
$\sum_{k=0}^n p_k=1$ and $p_n>0$, for some $n\ge 1$.
Then $F^{*\tau}\in\mathcal{L}$ if and only if $F^{*n}\in \mathcal{L}$.\\
(3) If, for some $n\ge 1$, ${\mathbf P}(\tau \ge n)>0$ and $F^{*n}\in\mathcal{L}$, and
if
\begin{equation}\label{1102}
\overline{G}(ax) = o\left(x^{1/2}\overline{F^{*n}}(x)\right)
\end{equation}
for any $a>0$,
then $F^{*\tau}\in \mathcal{L}$.\\
(4) Let $F$ and $L$ be two distributions
such that $F^{*2}\in\mathcal{L}$ and $\overline{F}(x)\sim \overline{L}(x)$.  Then $\overline{L^{*2}}(x)\sim
\overline{F^{*2}}(x)$ and, therefore, $L^{*2}\in\mathcal{L}$.
\end{Corol}

In order to illustrate the above results and to formulate the new ones,
we need further notion and notation.
Recall that a distribution $F$ is {\em dominatedly-varying-tailed}, denoted by $F\in\mathcal{D}$, if for some (or
equivalently, for all) $c\in(0, 1)$,
$$\limsup \overline{F}(cx)/\overline{F}(x)<\infty.$$
A distribution $F$ belongs to the {\em generalised long-tailed} distribution class $\mathcal{OL}$, if 
$\overline{F}(x)>0$ for all $x$ and if,
for any $c>0$,
$$C(F,c)=\limsup\overline {F}(x-c)/\overline {F}(x)<\infty.$$
The class $\mathcal{OL}$ is significantly broader than the class $\mathcal{L}$ and,
in particular, the class $\mathcal{OL}$ 
covers all classes $\mathcal{L} (\gamma )$, $\gamma \ge 0$.

The classes $\mathcal{D}$ and $\mathcal{OL}$ were introduced by \cite{F1969} and \cite{ShiWat2005}, respectively. Note that $\mathcal{OS}$ is a proper subclass of the class $\mathcal{OL}$, see e.g. \cite{ShiWat2005} or \cite{W2008}.

\begin{remark}\label{rem101}
Statement (1a) of Theorem \ref{th1} is quite general --
in particular, it may be applied in the case where $n=2$,
and $F_1$ is not long-tailed itself. We present two examples in the Appendix below. In Example \ref{exam501}, there are two distributions $F_1\in\mathcal{OL}$ and $F_2$ such that $\overline{F_2}(x)=o(\overline{F_1}(x))$; and in Example \ref{exam503}, there are two distributions $F_1\notin\mathcal{OL}$ and $F_2$ such that $\liminf\overline{F_2}(x)/\overline{F_1}(x)=0$ and $\limsup\overline{F_2}(x)/\overline{F_1}(x)=\infty$.
In both examples, $F_1\notin\mathcal{L}\cup\mathcal{D}$ and $F_1*F_2\in\mathcal{L}$.
\end{remark}

\begin{remark}\label{rem102}
In Corollary \ref{cor001}, parts (1) and (3), if $n\ge 2$, then distributions $F^{*k}$ may be not
long-tailed for $1\le k\le n-1$, in general -- see, e.g., families of distributions
$\mathcal{F}_i(0), i=1,2$ that are introduced below.
Therefore this result is a reasonable generalisation of Theorem 6 of
Leipus and \v{S}iaulys \cite{LeiSia2012}. Also,
Leipus and \v{S}iaulys \cite{LeiSia2012} require condition \eqref{1102}
with $n=1$ that is stronger than our condition if $F$ does not belong to the
class $\mathcal{OS}$.
\end{remark}

\begin{remark}\label{rem1020}
The results of part (1) of Theorem \ref{th1}
may be generalised onto the case of weakly dependent random variables.
Here is an example for $n = 2$, with a particular choice of a weak dependence structure of random variables. Let $X_i$ be a random variable with the distribution $F_i$ supported on whole real line, $i=1,2$. Assume that a random vector $(X_1,X_2)$ has  the two-dimensional Farlie-Gumbel-Morgenstern (FGM) joint distribution:
\begin{eqnarray}\label{503}
{\mathbf P}\Big(\bigcap_{i=1}^{2}\{X_i\le x_i\}\Big)=\prod_{i=1}^{2}F_i(x_i)(1+\theta_{12}\overline{F_1}(x_1)\overline{F_2}(x_2)),
\end{eqnarray}
where $\theta_{12}\ne 0$ is a constant such that $a=\mid\theta_{12}\mid\le1$.

For any $0<T_i\le\infty,i=1,2$, direct calculations show that
\begin{eqnarray}\label{504}
&&{\mathbf P}\Big(\bigcap_{i=1}^{2}\{X_i\in (x_i,x_i+T_i]\}\Big)=
\prod_{i=1}^{2}F_i(x_i,x_i+T_i]
\big(1+\theta_{12}\prod_{i=1}^{2}(1-\overline{F_i}(x_i+T_i)-\overline{F_i}(x_i))\big).
\end{eqnarray}
Then, by \eqref{504}, we have for $1\le i\neq j\le2$ and all $x_i,x_j$,
\begin{eqnarray}\label{505}
&&{\mathbf P}(X_i\in (x_i,x_i+T_i]|X_j=x_j)=F_i(x_i,x_i+T_i]\nonumber\\
&&\ \ \ \ \ \ \ \ \ \ \ \ \cdot\big(1+\theta_{12}(1-\overline{F_i}(x_i+T_i)-\overline{F_i}(x_i))(1-2\overline{F_j}(x_j))\big).
\end{eqnarray}

Take $a<1$. One may show that the 
statements (1a) and (1b) of Theorem \ref{th1}
still hold under new assumptions by simply following their proofs, with a suitable
use of
equalities \eqref{503} and \eqref{505}.
\end{remark}

\vspace{0.5cm}

Now we discuss the closedness property under convolution roots related to the class
$\mathcal{L}\cap\mathcal{OS}$.
We show that all three Conjectures 1-3 do not take place in the class $\mathcal{L}\cap\mathcal{OS}$. We provide precise examples and the intuition behind. All our examples involve
absolutely continuous distributions. In more detail, we introduce below two families of distributions,
$\mathcal{F}_1(0)$ and $\mathcal{F}_2(0)$, that have different properties and
are built up around random variables of the form
\begin{equation}\label{etau}
\xi = \eta (1+U)
\end{equation}
where $\eta$ has a discrete and heavy-tailed distribution and $U$ is an independent random
variable with a smooth distribution with bounded support. For simplicity, we assume $U$ to be uniformly distributed, but its
distribution may be taken from a larger class. Further,
classes $\mathcal{F}_1(0)$ and $\mathcal{F}_2(0)$ may be extended, thanks to part (4)
of Corollary \ref{cor001} on the tail-equivalence.

\begin{defin}\label{defin101}
Class $\mathcal{F}_1(0)$ is a 4-parametric family of distributions
$F=F(\alpha,b,t,A)$ of random variables
\be \label{exam201-1}
\xi=\eta(1+U^{1/b})^t
\ee
with density $f=f(\alpha,b,t,A)$. Here $\alpha \in [1/2,1)$, $b>0$ and $t\ge 1$  are constants.
Further, $\eta$ is a
discrete random variable with distribution ${\mathbf P}(\eta=a_n)=Ca_n^{-\alpha}$, where
$C=(\sum_{n=0}^{\infty}a_n^{-\alpha})^{-1}$ is the normalising constant and
a sequence $A=\{a_n\}$ is defined as follows. Let $r=1+1/\alpha >2$ and a constant $a>1$
be so large that $a^{r} > 2^{t+2}a$, then
$a_n =
a^{r^{n}}$ for $n= 0,1,\ldots.$
Finally, $U$ is a random variable having uniform distribution in the
interval $(0,1)$, and $U$ and $\eta$ are mutually independent.
\end{defin}


A number of ``good'' properties of the
class $\mathcal{F}_1(0)$ is given in the following theorem.
In particular, the theorem provides a negative answer
to Conjectures 1-3 related to the class $\mathcal{L}\cap\mathcal{OS}$.

\begin{thm}\label{th2} For any distribution $F\in \mathcal{F}_1(0)$, the following conclusions hold.\\
(1) $F$ is neither long-tailed nor generalised subexponential, while $F\in\mathcal{OL}$ and $F^{*n}\in \mathcal{L}\cap\mathcal{OS}\setminus \mathcal{S}$, for all $n\ge 2$.\\
(2) $F^{*\tau}\in \mathcal{L}\setminus \mathcal{S}$, for any counting random variable $\tau$ with distribution $G$ such that
${\mathbf P} (\tau \ge 2)>0$ and for any $a>0$
\be \label{thm101}
\overline{G}(ax) = o\big(x^{1/2}\overline{F^{*2}}(x)\big).
\ee
(2a) Further, if condition (\ref{thm101}) is replaced by the following: for any $0<\varepsilon<1$, there is an integer $M=M(\varepsilon)\ge2$ large enough such that
\be \label{thm102}
\sum_{n=M}^\infty p_n\overline{F^{*n}}(x)\le\varepsilon\overline{F^{*\tau}}(x),\ for\ all\ x\ge0,
\ee
then $F^{*\tau}\in \mathcal{L}\cap\mathcal{OS}$.\\
(2b) Assume now that ${\mathbf E}\tau<\infty$. Then \eqref{thm102} implies that
\be \label{thm1020}
\liminf\overline{F^{*\tau}}(x)/\overline{F}(x)
=
{\mathbf E}\tau
\ge
2\liminf\overline{F^{*\tau}}(x)/\overline{F^{*2}}(x)
\ge
2\sum_{m=1}^\infty(p_{2m}+p_{2m+1})m.
\ee
Further, if condition
\begin{equation}\label{rem1040}
\sum_{m=1}^\infty\Big(\sum_{k=2(m-1)+1}^{2m}p_k\Big)\big(C^*(F^{*2})-1+\varepsilon_0\big)^m<\infty
\end{equation}
holds 
for some $\varepsilon_0>0$, then
\be \label{thm1021}
2\limsup\overline{F^{*\tau}}(x)/\overline{F^{*2}}(x)\le
2\sum_{m=1}^\infty m(p_{2m-1}+p_{2m})(C^*(F^{*2})-1)^{m-1}
<\infty
\ee
while $\limsup\overline{F^{*\tau}}(x)/\overline{F}(x)=\infty$.

(3) For any distribution $F\in\mathcal{F}_1(0)$, there is an infinitely divisible distribution $H$ such that $F$ is generated by its L$\acute{e}$vy measure and the following holds:  $H\in
(\mathcal{L}\cap\mathcal{OS})\setminus\mathcal{S}$, while $F$ is neither
long-tailed nor generalised subexponential.

\end{thm}

\begin{remark}
Assume a random variable $\tau$ has a Poisson distribution with parameter
$\mu ={\mathbf E} \tau$. Let $r=(C^*(F^{*2})-1)^{1/2}$.
Direct computations show
that the lower bound in \eqref{thm1020} is equal to
$$
\mu + (1-e^{-2\mu})/2,
$$
and the upper bound in  \eqref{thm1021} is equal to
$$
\frac{\mu+1}{2r}\left( e^{\mu (r-1)}-e^{-\mu (r+1)}\right)
+\frac{\mu}{2} \left(e^{\mu (r-1)} + e^{-\mu (r+1)}\right).
$$
The same (lower and upper)
bounds hold for  the lower and upper limits of $2\overline{H}(x)/\overline{F^{*2}}(x)$
in
part (3) of Theorem \ref{th2}. In this case, $\mu$ is precisely given in the proof, see Section 3.
\end{remark}

\begin{lemma}\label{gengen}
The following
condition
implies \eqref{thm1020}: there exist $n\ge 1$ and $\varepsilon_0>0$ such that
\begin{equation}\label{rem104}
\sum_{m=1}^\infty\Big(\sum_{k=(m-1)n+1}^{mn}p_k\Big)\big(C^*(F^{*n})-1+\varepsilon_0\big)^m<\infty.
\end{equation}
\end{lemma}
One can see that condition \eqref{rem1040} is a particular case of
condition \eqref{rem104}, with $n=2$.

\begin{remark}\label{rem103}
Condition \eqref{rem104} holds if a distribution $G$ is either
Poisson ($p_k=\lambda^k e^{-\lambda}/(k!),\ k=0,1,\ldots$)
or Geometric
($p_k=qp^k,\ k=0,1,\ldots$, with $p <1/(C^*(F^{*n})-1+\varepsilon_0)$, for
some $\varepsilon_0>0$).
Note that \eqref{rem104} is a natural generalisation of
the classical sufficient condition for subexponentiality of a random sum
(where $n=1$ and $C^*(F)=2$), see e.g.
Theorem 4 in \cite{ChoNeyWai1973a} or Theorem 3 and its Remark in \cite{EmbGolVer1979}.
Clearly, a distribution $G$ satisfying \eqref{rem104} is light-tailed.
\end{remark}

Here is an example
of
a heavy-tailed distribution $G$ that satisfies condition \eqref{thm102}.
\begin{exam}\label{exam505} Let $n=1$. Assume $F\in\mathcal{D}$,
then by Theorem 3 of Daley et al. \cite{DOV2007}, there are two positive constants $C$ and $\alpha$ such that
\begin{eqnarray}\label{509}
\sup_{x\ge0}\overline{F^{*k}}(x)/\overline{F}(x)\le Ck^{\alpha},\ for\ all\ k\ge1.
\end{eqnarray}
Take a counting random variable $\tau$ with distribution $G$ given by
${\mathbf P} (\tau=k)=
p_k=Kk^{-\beta}$ for some $\beta>\alpha+2$, where
$K=\left(\sum_{k=1}^{\infty} k^{-\beta}\right)^{-1}$ is the normalising constant.
Clearly, condition (\ref{thm102}) takes place and $G$ is a heavy-tailed distribution. However, condition (\ref{rem104}) in Remark \ref{rem103}
does not hold.
\end{exam}

\begin{remark}\label{rem1030}
Note that all distributions $F$ considered in Theorem \ref{th2} are generalised
long-tailed,
that is
$\mathcal{F}_1(0)\subset\mathcal{OL}$. One may guess that such a
condition may be essential for $F^{*2}$ to be long-tailed. However, this is not
the case: we introduce below another family $\mathcal{F}_2(0)$ of heavy-tailed distributions $F$ such that $F\notin\mathcal{OL}$ while $F^{*2}\in\mathcal{L}$
and, moreover, $F^{*2}\in\mathcal{OS}$.
\end{remark}

\begin{defin}\label{defin102}
Class $\mathcal{F}_2(0)$ is a 3-parametric family of heavy-tailed distributions
$F=F(\alpha,t,A)$ of random variables
\be \label{def201-1}
\xi=\eta^{1/t}(1+U)^{1/t}
\ee
with density $f=f(\alpha,t,A)$.
Here $t\in (1,2)$, $\alpha \in ((1-t)/t,1/t)$ and the sequence $A=\{a_n\}$ and
random variables $\eta$ and $U$ are defined
as in Definition \ref{defin101}.
\end{defin}

Properties of the class $\mathcal{F}_2(0)$ are summarised in the following proposition.

\begin{pron}\label{classF2}
Let $F\in \mathcal{F}_2(0)$, then
$F^{*n}\in \mathcal{L}\setminus \mathcal{S}$,
for all $n\ge 2$.
Further, for any $n\ge 2$, $F^{*n}\in\mathcal{OS}$ when $\alpha\in[1/2,1/t)$ and
$F^{*n}\notin\mathcal{OS}$ when $\alpha\in((t-1)/t,1/2)$, while $F\notin\mathcal{OL}$,
and therefore $F\notin\mathcal{L}\cup\mathcal{D}$.
\end{pron}

\begin{remark}
In addition, for the class $\mathcal{F}_2(0)$ with $\alpha\in[1/2,1/t)$, the natural analogues of statements (2) and (3) of Theorem \ref{th2} do hold.
\end{remark}

The {\bf proof of Proposition \ref{classF2}} is quite similar to that of Theorem \ref{th2}.
For the sake of completeness, we decided to give it in Subsection 4.2 of Appendix.

Theorem \ref{th2} and Proposition \ref{classF2} provide a good number of new examples of distributions
from the classes ${\mathcal L}\setminus {\mathcal S}$
and  $(\mathcal {L}\cap\mathcal{OS})\setminus {\mathcal S}$.
\vspace{0.2cm}

\begin{remark}\label{Watanabe}
Watanabe and Yamamuro \cite{WY2010} commented in
Remark 2.3 that Shimura and Watanabe \cite{ShiWat2005b}
provided a counter-example to Conjecture 1.
Also, \cite{WY2010} pointed out that \cite{ShiWat2005b} did not find an answer to the corresponding Conjectures 2-3 related to
distributions of random sums (compound distribution or random convolution) and infinitely divisible distribution.
In addition, \cite{WY2010} stated that the class $\mathcal{OS}$ is not closed under
convolution roots, but we did not find any corresponding result for the intersection of the classes $\mathcal{L}\cap\mathcal{OS}$.

Recently we were in touch with Dr Shimura who has sent
us privately an unpublished English translation of Research
Report \cite{ShiWat2005b}. We have found that the counter-example there
seems to be correct,  but is described implicitly, so it is difficult to follow. Also, the example
relates to a distribution that is neither absolutely continuous
nor discrete.
\end{remark}

Finally, we show that
the long-tailedness property is preserved under
convolution roots within the class $\mathcal{OS}$. Namely, the following result holds.

\begin{pron}\label{pron1}
If $F\in\mathcal{OS}$, then $F\in\mathcal{L}$ if and only if $F^{*2}\in\mathcal{L}$.
\end{pron}

\section{Proofs of Theorem \ref{th1} and Corollary \ref{cor001}}
\setcounter{thm}{0}\setcounter{Corol}{0}\setcounter{lemma}{0}\setcounter{pron}{0}
\setcounter{remark}{0}\setcounter{exam}{0}\setcounter{equation}{0}\setcounter{property}{0}

In order to prove Theorem \ref{th1}, we first recall a number of known properties of long-tailed distributions. We consider here distributions on the whole real line.

The definition of the class $\mathcal{L}$ and the diagonal argument lead to the following result.

\begin{property}\label{property201}
Distribution $F$ is long-tailed if and only if there exists
a monotone increasing function $h(x)\uparrow\infty$ such that $h(x)<x$ and
$F(x-h(x),x+h(x)] = o(\overline{F}(x))$ (then we say that $\overline{F}$ is
$h$-{\it insensitive}).
\end{property}

See, e.g., \cite{FosKorZac2013}, Chapter 2 for Property 1 and for
$h$-insensitivity and other properties of class $\mathcal{L}$. Further, \cite{EmbGol1980} and \cite{EmbGol1982} show that the class $\mathcal{L}$ is closed under convolution and mixture.

\begin{property}\label{property202}
Let $F_1$ and $F_2$ be two distributions.\\
(1) Assume $F_1\in\mathcal{L}$. Then $F_1*F_2\in \mathcal{L}$
if either (a) $F_2\in\mathcal{L}$ or (b) $\overline{F_2}(x)=o(\overline{F_1}(x))$.
In the latter case, $\overline{F}_1(x) \sim \overline{F_1*F_2}(x)$.\\
(2) If $F_1,F_2\in\mathcal{L}$, then $pF_1+(1-p)F_2\in\mathcal{L}$, for any $p\in [0,1]$.
\end{property}

Albin \cite{Alb2008} and then Leipus and \v{S}iaulys \cite{LeiSia2012} extended
Property 2 (1) onto random convolutions.

\begin{property}\label{property203}
If $F\in\mathcal{L}$ and if
\eqref{1102}
holds for $n=1$ and for all $a>0$,
then $F^\tau\in \mathcal{L}$.
\end{property}

We proceed now with the
\noindent{\bf Proof of Theorem \ref{th1}.}

{\bf Proof of (1a).} First, we prove (\ref{1011}) for $i=1$.
By $H_n\in\mathcal{L}$,
we may choose $h(x)\uparrow\infty$ such that $\overline{H}_n$ is
$nh$-insensitive.
Then, by Property 1,
$$
{\mathbf P} \big(S_n\in (x-nh(x), x+nh(x)]\big) = o\big(\overline{H_n}(x)\big).
$$
Note that
$$
{\mathbf P}\big (S_n\in (x-nh(x), x+nh(x)]\big)\ge
{\mathbf P} \big(X_1\in (x-h(x), x+h(x)]\big)\cdot \prod_{j=2}^n
{\mathbf P} \big(X_j \in (-h(x), h(x)]\big)
$$
and
$${\mathbf P} \big(-h(x)< X_j \le h(x)\big)\to 1, \quad j=2,\ldots,n.$$
Then the first part of \eqref{1011} follows.
Since $H_1=F_1$, the second part follows too.

If $i>1$, then the proof of the first part of \eqref{1011} is the same.
For the second part, we may represent $S_n$ as a sum of mutually independent random
variables $S_n = S_i+X_{i+1}+\ldots + X_n$ and apply the arguments from
above.
\hfill$\Box$

{\bf Proof of (1b).} It is enough to prove the result for $m=1$ and $n=2$, and then use the
induction argument.
First, by monotonicity of distribution functions and since $F_1$ is
long-tailed, we
may obtain that  $F_2(x-c,x+c] = o(\overline{F_1}(x))$ for any
$c>0$ and, therefore,
$$
\alpha_c (x)=: \sup_{y\ge x}\Big(F_2(x-c,x+c]/\overline{F_1}(y)\Big) \downarrow 0.
$$
Then one can use the diagonal argument to conclude that there
exists a positive function $h_1(x)\uparrow\infty$ such that
\begin{equation}\label{2h1}
F_2(x-2h_1(x),x+2h_1(x)]
 = o(\overline{F_1}(x)).
\end{equation}
Further, since $F_1$
is long-tailed, one can find a function $h_2(x)\uparrow\infty$ such that $\overline{F_1}$
is $2h_2$-insensitive. Let $h(x)=\min (h_1(x),h_2(x))$. Then $\overline{F_1}$
is $2h$-insensitive and \eqref{2h1} holds with $h$ in place of $h_1$.

Let $X_1,X_2$ be two independent random variables where $X_1$ has distribution $F_1$ and
$X_2$ has distribution $F_2$. Then, for any $c>0$ and for $x$ such that $h(x)>c$,
\begin{eqnarray*}
&&F_1*F_2 (x-c,x+c] =
{\mathbf P} (X_1+X_2 \in   (x-c, x+c]) \\
&\le &
F_2(x-2h(x), x] + \Big(
\int_{-\infty}^{x-2h(x)} + \int_{x}^{\infty}
\Big)
F_2(dy)F_1
(x-y-c, x-y+c].
\end{eqnarray*}

There are three terms on the right-hand side.
The first term  is $o(\overline{F_1}(x))$, by condition \eqref{2h1}.
It is also $o(\overline{F_1*F_2}(x))$ since
$\overline{F_1*F_2}(x) \ge \overline{F_1}(x_0)\overline{F_¬ß}(x-x_0) \sim
\overline{F_1}(x_0)\overline{F_1}(x)$, where $x_0$ is any number such that
$\overline{F_2}(x_0)>0$.

Then the second term is not bigger than
$$
\alpha_c (2h(x)-c)  \int_{-\infty}^{x-2h(x)}F_2(dy) \overline{F_1
(x-y)} \le   \alpha_c (2h(x)-c) \overline{F_1*F_2}(x) = o(\overline{F_1*F_2}(x)).
$$

Finally, the last term is not bigger than
\begin{eqnarray*}
&& \sum_{k=0}^{\infty} F_2(x+kc, x+(k+1)c]
F_1(-(k+2)c, -(k-1)c] \\
&\le & 3 \sup_{y\ge x} F_2(y,y+c]F_1(c)\\
&=& o(\overline{F_1}(x))=o(\overline{F_1*F_2}(x)).
\end{eqnarray*}

Thus $F_1*F_2\in\mathcal{L}$. \hfill$\Box$

{\bf Proof of (2).}
Assume first that $H_n\in\mathcal{L}$. Then, by property (1a),
${H_k}(x-c,x+c]=o(\overline{H_n}(x))$,
for all $k=1,\ldots,n$ and for any fixed $c>0$. Then
$$ H_{\tau}(x-c,x+c]
 =
\sum_{k=1}^n p_k H_k(x-c,x+c]
= o(\overline{H_n}(x)),
$$
and $H_{\tau}\in \mathcal{L}$ follows.

Vice versa, if $H_{\tau}\in\mathcal{L}$, then
$$H_k(x-c,x+c]
\le {H}_{\tau}(x-c,x+c])/p_k =o(\overline{H}_{\tau}(x))
$$
for each $k$ such that $p_k>0$ and, in particular, for $k=n$. Let $x_1,\ldots,x_n$ be positive numbers
such that $\overline{F}_i(x_i)>0$.
Clearly,
$\overline{H}_{n}(x)\ge\overline{H_k}(x-\sum_{i=k+1}^n x_i)\prod_{i=k+1}^n\overline{F_i}(x_i),$
for all $k=1,\ldots,n-1$ and then
$$
\overline{H}_{\tau}(x) \le \sum_{k=1}^n p_k \overline{H_k}\Big(x-\sum_{i=k+1}^n x_i\Big)
\le
\overline{H_{n}}(x)/\Big(\prod_{i=1}^n\overline{F_i}(x_i)\Big).$$
Thus $H_n\in\mathcal{L}$ follows from
\begin{eqnarray*}
H_n (x-c,x+c]
=o(\overline{H}_{\tau}(x))= o(\overline{H_n}(x)).
\end{eqnarray*}
\hfill$\Box$

{\bf Proof of (3).} We may assume, without loss of generality, that $p_n={\mathbf P} (\tau =n)>0$.
Further, we may assume that ${\mathbf P}(\tau >n)>0$ -- otherwise
the result follows from the previous statement.

Let $P_n= {\mathbf P} (\tau \le n)=\sum_{k=0}^np_k$ and
$Q_n={\mathbf P} (\tau >n)=\sum_{k=n+1}^{\infty}p_k$.
Further, let $H^{(1)}(x)=\sum_{k=1}^n p_k H_k(x)/P_n$
and $H^{(2)}(x)=\sum_{k=n+1}^{\infty} p_k H_k(x)/Q_n$. Since $H=P_nH^{(1)}+Q_nH^{(2)}$,
it is enough to show that both $H^{(1)}$ and $H^{(2)}$ are long-tailed -- see Property 2 (2).
By the previous statement (2), we have $H^{(1)}\in\mathcal{L}$. Then the argument from
\cite{LeiSia2012} implies
\begin{eqnarray*}
\overline{H^{(2)}}(x-1)
&=& \sum_{n+1\le k \le (x-x_0)/(x_0-1)} p_k\overline{H_{k}}(x-1)/Q_n
+\sum_{k > (x-x_0)/(x_0-1)} p_k\overline{H_{k}}(x-1)/Q_n\\
&\le & (1+\varepsilon )\overline{H^{(2)}}(x) +
\sum_{k > (x-x_0)/(x_0-1)} p_k{H_{k}}((x-1,x])/Q_n\\
&\le &
(1+\varepsilon )\overline{H^{(2)}}(x) +
\frac{\overline{G}((x-x_0)/(x_0-1))}{Q_n\sqrt{(x-x_0)/(x_0-1)}} \\
&=&
(1+\varepsilon )\overline{H^{(2)}}(x) +
o(\overline{H^{(2)}}(x)).
\end{eqnarray*}
Since $\varepsilon >0$ is arbitrary, the distribution $H^{(2)}$ is long-tailed.\hfill$\Box$

{\bf Proof of (4)}.
By part (1) of the theorem, there exists a function $h(x)\uparrow \infty$ such that
$\overline{F_1*F_2}$ is $h$-insensitive and
$\overline{F_1}(x-h(x))-\overline{F_1}(x)=o(\overline{F_1*F_2}(x))$.
Then
$$
\int_{x-h(x)}^x F_1(dy) \overline{F}_2(x-y)
\le F_1(x-h(x),x]
=o(\overline{F_1*F_2}(x))
$$
and, similarly, $\int_{x-h(x)}^x F_1(dy) \overline{L_2}(x-y)
= o(\overline{F_1*F_2}(x))$.\\

Next,
$$
 \overline{F_1*F_2}(x)/\overline{F}_1(x)\sim
 \overline{F_1*F_2}(x-h(x))/\overline{F}_1(x)
\ge \overline{F}_1(x)\overline{F}_2(-h(x)) /\overline{F}_1(x) = 1-o(1)
$$
and then
\begin{eqnarray*}
\overline{F}_1(x) &\ge& \left(
\int_x^{x+h(x)} + \int_{x+h(x)}^{\infty} \right)
 F_1(dy) \overline{F}_2(x-y)\\
&\ge& o(\overline{F_1*F_2}(x)) +\int_{x+h(x)}^{\infty}F_1(dy)\overline{F}_2(-h(x))\\
&=& o(\overline{F_1*F_2}(x)) + (\overline{F}_1(x+h(x))
-\overline{F}_1(x)) (1+o(1))+\overline{F}_1(x)(1+o(1))
\\
&=& o(\overline{F_1*F_2}(x))+ \overline{F}_1(x).
\end{eqnarray*}
Therefore, $\int_x^{\infty} F_1(dy)\overline{F}_2(x-y) =
\overline{F}_1(x) + o(\overline{F_1*F_2}(x))$ and the same holds with
$L_2$ in place of $F_2$ in the left-hand side of the latter
equality.

Further, due to the monotonicity of distribution functions,
$\overline{F}_2(x-y)\sim \overline{L_2}(x-y)$ uniformly in
${x-y}\ge h(x)$. Therefore
$$
\int_{-\infty}^{x-h(x)}F_1(dy)\overline{F}_2(x-y) \sim
\int_{-\infty}^{x-h(x)}F_1(dy)\overline{L_2}(x-y).
$$

Finally,
\begin{eqnarray*}
\overline{F_1*L_2}(x) &=&
\int_x^{\infty} F_1(dy)\overline{L_2}(x-y) +\int_{-\infty}^{x-h(x)} F_1(dy)
\overline{L_2}(x-y) +
\int_{x-h(x)}^x F_1(dy) \overline{L_2}(x-y)\\
&\sim&
\int_x^{\infty} F_1(dy)\overline{F}_2(x-y)
 + \int_{-\infty}^{x-h(x)}F_1(dy) \overline{F}_2(x-y)
+o(\overline{F_1*F_2}(x))\\
&\sim &
\overline{F_1*F_2}(x),
\end{eqnarray*}
and therefore $F_1*L_2\in\mathcal{L}$.\hfill$\Box$

\vspace{0.2cm}

\noindent{\bf Proof of Corollary \ref{cor001}.}
We need to prove statements (1), (3) and (4) only.

{\bf Proof of (1).} If $F^{*2}\in\mathcal{L}$, then by statement (1a) of Theorem \ref{th1}, $\overline{F}(x-t)-\overline{F}(x+t)
= o(\overline{F^{*2}}(x))$ for any $t>0$. Further, by statement (1b) of
Theorem \ref{th1}, we have $F^{*3}= F^{*2}*F \in \mathcal{L}$. Then Property 2 and
the induction argument complete the proof. \hfill$\Box$\\

{\bf Proof of (3).} Condition \eqref{concentration} follows from \cite{Pet1995}, Theorem 2.22.
So we have to verify \eqref{tail} only.
Due to Lemma 2.1 from \cite{Alb2008} or Lemma 4 from \cite{LeiSia2012}, for any long-tailed distribution $V$
and for any $\varepsilon >0$,
there is
$x_0>1$ such that, for all
$i\ge1$,
\begin{eqnarray}\label{304}
\sup_{x\ge
n(x_0-1)+x_0}\overline{V^{*n}}(x-1)/\overline{V^{*n}}(x)\le
1+\varepsilon.
\end{eqnarray}
Clearly, if there are, say, $m$ long-tailed distributions $V_1,\ldots,V_m$, then
\eqref{304} holds again for some $x_0>1$ and for any $V_i$ in place of $V$. Using similar
arguments, one can also show that , for any $i\ge 1$, inequalities \eqref{304}
hold for $U_n$ in place of $V^{*n}$ where $U_n$ is any convolution of $n$ distribution
functions taken from the set $\{V_1,\ldots, V_m\}$ -- namely, $U_n=V_1^{*j_1}*\ldots
*V_m^{*j_m}$ where $j_1+\ldots +j_m=n$. As the corollary, we may take $m=n$ and then
$V_l=F^{*(n+l)}$, for $l=1,\ldots,n$, to conclude that inequalities \eqref{304}
continue to hold for $i\ge n$, with $F^{*n}$ in place of $V^{*n}$.\hfill$\Box$
\vspace{0.2cm}

{\bf Proof of (4).} We have to apply part (4) of Theorem \ref{th1} twice, first to move
from
$F^{*2}$ to $F*L$ and then from $F*L$ to $L^{*2}$.
\hfill$\Box$

\section{Proofs of Theorem \ref{th2}, Proposition \ref{pron1} and Lemma \ref{gengen}}
\setcounter{thm}{0}\setcounter{Corol}{0}\setcounter{lemma}{0}\setcounter{pron}{0}
\setcounter{remark}{0}\setcounter{exam}{0}\setcounter{equation}{0}\setcounter{property}{0}

We start with a simple auxiliary result.
\begin{lemma}\label{lem301}
Assume that a distribution $F$ is absolutely continuous with density $f$. If
\begin{equation}\label{301}
f(x) = o(\overline{F}(x))\ a.e.,
\end{equation}
then $F$ is long-tailed.\hfill$\Box$
\end{lemma}

\noindent{\bf Proof.} Indeed, let $\varepsilon(x) = \sup_{y\ge x} f(y)/
\overline F(y)$. Since $\varepsilon(x)\downarrow 0$, we have
$$\overline F(x+1)\leq\overline F(x)=\overline {F}(x+1)+\int_{x}^{x+1}f(y)dy
\leq\overline F(x+1)+\varepsilon(x)\overline F(x),$$
and the result follows. \hfill$\Box$
\vspace{0.2cm}

\noindent{\bf Proof of Theorem \ref{th2}.} Start with {\bf Proof of (1).} Recall that $F\notin\mathcal{L}$ implies
with necessity that $F^{*n}\notin\mathcal{S}$ for all $n\ge2$.
Then, by Corollary \ref{cor001} of the present paper and Proposition 2.6 from [25], we only need to prove
that $F\notin\mathcal{L}\cup\mathcal{OS}$, $F\in\mathcal{OL}$ and $F^{*2}\in\mathcal{L}\cap\mathcal{OS}$.

First, we find closed-form representations for distribution $F$ and its density $f$.
Clearly, $\eta\leq\xi\leq2^{t}\eta$. Since
\be
\label{exam201-2}\textbf{E}\eta^{s}=C\sum_{n=0}^{\infty}a_n^{s-\alpha}<\infty
\ee
if and only if $s<\alpha$, the same holds for $\xi$, and distribution $F$ is
heavy-tailed with infinite mean. Further, by \eqref{exam201-1} we have,
for $n\ge1$,
\begin{eqnarray*}
&&F(a_n,x]\textbf{I}(x\in [a_n,a_{n+1}))\\
&=&
\textbf{P}(\eta=a_n)\textbf{P}(1<(1+U^{1/b})^t\le x/a_n)
\Big(\textbf{I}(x\in[a_n,2^ta_n))+\textbf{I}(x\in[2^ta_n,a_{n+1}))\Big)\\
&=&
Ca_n^{-\alpha}\Big((a_n^{-1}x)^{1/t}-1\Big)^{b}\textbf{I}(x\in[a_n,2^ta_n))
+C a_n^{-\alpha}\textbf{I}(x\in[2^ta_n,a_{n+1})).
\end{eqnarray*}
Then
\begin{eqnarray}\label{exam201-4}
f(x)=Cbt^{-1}\sum_{n=0}^{\infty}x^{1/t-1}a_n^{-\alpha-1/t}
\Big((xa_n^{-1})^{1/t}-1\Big)^{b-1}\textbf{I}(x\in [a_n,2^ta_n))
\end{eqnarray}
and
\begin{eqnarray}\label{exam201-5}
\overline{F}(x)&=&\textbf{I}(x<
a_0)+\sum\limits_{n=0}^{\infty}\Big(\textbf{P}(\xi\in(a_n,a_{n+1}])-\textbf{P}(\xi\in(a_n,x])
+\textbf{P}(\xi>a_{n+1})\Big)\textbf{I}(x\in[a_n,2^ta_n))\nonumber\\
& &\ \ \ \ \ \
+\sum\limits_{n=0}^{\infty}\Big(\textbf{P}(\xi\in(2^ta_n,a_{n+1}])-\textbf{P}(\xi\in(2^ta_n,x])
+\textbf{P}(\xi>a_{n+1})\Big)\textbf{I}(x\in[2^ta_n,a_{n+1}))\nonumber\\
&=&\textbf{I}(x<
a_0)+\sum\limits_{n=0}^{\infty}\Bigg(\bigg(\sum\limits_{i=n}^{\infty}
Ca_i^{-\alpha}-Ca_n^{-\alpha}\Big((x/a_n)^{1/t}-1\Big)^{b}\bigg)
\textbf{I}(x\in[a_n,2^ta_n))\nonumber\\
& &\ \ \ \ \ \
+\sum\limits_{i=n+1}^{\infty}Ca_i^{-\alpha}\textbf{I}(x\in[2^ta_n,
a_{n+1}))\Bigg),\ \ \ x\in(-\infty,\infty).
\end{eqnarray}

Now, we prove that $F\in\mathcal{OL}\backslash\mathcal{L}$. Note that
$a_{n+1}a_n^{-2}\to\infty$ as $n\rightarrow\infty$, so for any $K>0$,
\be\label{exam201-6}\sum_{n\geq N}a_n^{-K}\sim a_N^{-K},
\ \ \overline F(a_n)\sim \textbf{P}(\eta=a_n)\ \quad \mbox{and} \quad
\textbf{P}(\eta>a_n)=o(\textbf{P}^2(\eta=a_n)).\ee
From (\ref{exam201-5}) and (\ref{exam201-6}), we have
\begin{eqnarray*}
\overline F(2^ta_n)\sim
\textbf{P}(\eta=a_{n+1})=Ca_{n+1}^{-\alpha}=Ca_n^{-1-\alpha}
\end{eqnarray*}
and
\begin{eqnarray*}
\overline F(2^ta_n-1)-\overline F(2^ta_n)=
Ca_n^{-\alpha}\lo(1-\lo((2^t-a_n^{-1})^{1/t}-1\ro)^b\ro)
\sim Cbt^{-1}2^{-t+1}a_n^{-\alpha-1},
\end{eqnarray*}
as $n\to\infty$. Therefore
\be\label{exam201-7}
\limsup_{x\to\infty} \overline F(x-1)/\overline F(x)=bt^{-1}2^{-t+1}+1,
\ee
so $F\notin\mathcal{L}$, but $F\in\mathcal{OL}$.

Next, we prove that $F^{*2}\in\mathcal{L}$.
Let $(\eta_i,U_i)$, $i = 1,2$ be two independent
copies of $(\eta,U)$, and let $\xi_i = \eta_i(1+{U_i}^{1/b})^t$ and $S_2= \xi_1
+ \xi_2$. The random variable $S_2$ has an absolutely continuous
distribution, say, $H = F^{*2}$ with density function
\begin{eqnarray}\label{exam201-8}
h(x)=\int_{0}^{x}f(y)f(x-y)dy=2\int_{x/2}^{x}f(y)f(x-y)dy,\
x\in(-\infty,\infty).
\end{eqnarray}
Clearly, $h(x) > 0$ if and only if $a_n+a_0 < x < 2^{t+1}a_n$, for $n =
0,1,\ldots$.
According to Lemma \ref{lem301}, it is enough to show that
\be
\label{exam201-9}h(x)=o(\overline H(x)).
\ee
We consider two cases: (i) $
x\in J_{n,1}=[a_n+a_0,3\cdot2^{t-1}a_n)$ and (ii) $ x\in
J_{n,2}=[3\cdot2^{t-1}a_n,2^{t+1}a_n)$
for $n = 0,1,\ldots$.

In the case (i), representations (\ref{exam201-4}) and (\ref{exam201-8})
lead to
\begin{eqnarray*}
h(x)&\le&2Cbt^{-1}a_n^{-\alpha-1/t}\int_{a_n}^{2^{t}a_n}y^{1/t-1}
\Big((ya_n^{-1})^{1/t}-1\Big)^{b-1}f(x-y)dy
\\
&\le&2Cbt^{-1}a_n^{-\alpha-1}\int_{a_n}^{2^{t}a_n}f(x-y)dy\le2Cbt^{-1}a_n^{-\alpha-1},
\end{eqnarray*}
while by (\ref{exam201-5})
\begin{eqnarray*}
\overline H(x)\ge\overline F^2(x/2)\ge\overline F^2(3\cdot2^{t-2}a_n)\ge
C^2a_n^{-2\alpha}\bigg(1-\Big(2\cdot(3\cdot4^{-1})^{1/t}-1\Big)^b\bigg)^2.
\end{eqnarray*}
Since $\alpha < 1$, $\sup_{x\in J_{n,1}}h(x)/H(x) \to 0$ as $n\to\infty$.

In the case (ii), representations (\ref{exam201-4}) and (\ref{exam201-8})
imply that
\begin{eqnarray*}
h(x)&=&2Cbt^{-1}a_n^{-\alpha-1/t}\int_{2^{-1}x}^{2^ta_n}
y^{1/t-1}\Big((ya_n^{-1})^{1/t}-1\Big)^{b-1}f(x-y)dy
\\
&\le&2Cbt^{-1}a_n^{-\alpha-1}\int_{x-2^ta_n}^{x/2}f(y)dy \\
&\le&2Cbt^{-1}a_n^{-\alpha-1}\overline
F(2^{t-1}a_n)\le2C^2bt^{-1}a_n^{-2\alpha-1},
\end{eqnarray*}
and by (\ref{exam201-5}) we get that
\begin{eqnarray*}
\overline H(x)\ge\overline F(x)\ge Ca_n^{-\alpha-1}.
\end{eqnarray*}
Then again $\sup_{x\in J_{n,2}}h(x)/H(x) \to 0$ as $n\to\infty$.

We may conclude that $(\ref{exam201-9})$ holds, therefore
$F^{*2}\in\mathcal{L}$.

In order to prove $F^{*2}\in\mathcal{OS}$, we only need to show that
\be\label{309}
T(x)=\int_{x/2}^{x}\overline{H}(x-y)h(y)dy=O(\overline H(x)).
\ee
It is clear that $T(x)>0$ if and only if $a_n+a_0 < x < 2^{t+2}a_n$, for $n =
0,1,\cdot\cdot\cdot$.

By (\ref{exam201-4}) and (\ref{exam201-8}), it is easy to see that,
for $n=0,1,\ldots$ , if $x\in [a_n+a_0,2^{t}a_n)$, then
\begin{eqnarray}\label{exam501-2}
h(x)=
2\int_{a_n}^{x}f(x-y)f(y)dy\le2Cbt^{-1}a_n^{-\alpha-1},
\end{eqnarray}
and if $x\in [2^{t}a_n,2^{t+1}a_n)$, then
\begin{eqnarray}\label{exam501-3}
h(x)=
2\int_{x/2}^{2^{t}a_n}f(x-y)f(y)dy\leq2Cbt^{-1}a_n^{-\alpha-1}\overline F(x-2^{t}a_n).
\end{eqnarray}
Then we estimate $T(x)$ separately in three cases: (i) $
x\in [a_n+a_0,3\cdot2^{t-1}a_n)$, (ii) $ x\in [3\cdot2^{t-1}a_n,2^{t+1}a_n)$ and (iii) $ x\in [2^{t+1}a_n,2^{t+2}a_n)$ for $n = 0,1,\ldots$.

In the case (i), representations (\ref{exam201-5}), (\ref{exam501-2}) and (\ref{exam501-3})
lead to
\begin{eqnarray*}
T(x)/\overline{H}(x)&\le&\max_{y\in[x/2,x]}\{h(y)\}\int_{x/2}^{x}\overline{H}(x-y)dy/
\overline{F}^2(2^{-1}x)\nonumber\\
&\le&2Cbt^{-1}a_n^{-\alpha-1}\int_{0}^{3\cdot2^{t-2}a_n}\overline{H}(y)dy/\overline{F}^2(3\cdot2^{t-2}a_n)
<\infty.
\end{eqnarray*}

In the case (ii), representations (\ref{exam201-5}), (\ref{exam501-2}) and (\ref{exam501-3}) imply that
\begin{eqnarray*}
T(x)/\overline{H}(x)&\le&\Big(\int_{x/2}^{2^{t}a_n}
+\int_{2^{t}a_n}^{x}\Big)\overline{H}(x-y)h(y)dy/
\Big(\overline{F}(x)+\int_{x/2}^{x}\overline{F}(x-y)F(dy)\Big)\nonumber\\
&\lesssim&2bt^{-1}\int_{x/2}^{2^{t}a_n}\overline{H}(x-y)dy/
\Big(1+\int_{x/2}^{2^{t}a_n}\overline{F}(x-y)dy\Big)\nonumber\\
& &\ \ \ \ \ \ +2bt^{-1}\int_{2^{t}a_n}^{x}\overline{H}(x-y)\overline{F}(y-2^{t}a_n)dy\nonumber\\
&\leq&4bt^{-1}\int_{x-2^{t}a_n}^{x/2}\Big(\overline{F}(y)+\int_{y/2}^{y}\overline{F}(y-z)F(dz)\Big)dy/
\Big(1+\int_{x-2^{t}a_n}^{x/2}\overline{F}(y)dy\Big)\nonumber\\
& &\ \ \ \ \ \ +2bt^{-1}\Big(\int_{0}^{x/2-2^{t-1}a_n}+
\int_{x/2-2^{t-1}a_n}^{x-2^{t}a_n}\Big)\overline{H}(y)\overline{F}(x-2^{t}a_n-y)dy\nonumber\\
&\leq&4bt^{-1}\Big(1+2Cbt^{-1}a_n^{-\alpha-1}(2^{t}a_n-2^{-1}x)
\int_{0}^{4^{-1}x}\overline{F}(z)dz\Big)
\nonumber\\& &\ \ \ \ \ \ +2bt^{-1}\Big(\overline{F}(2^{t-2}a_n)\int_{0}^{2^{t-1}a_n}\overline{H}(y)dy+
\overline{H}(2^{t-2}a_n)\int_{0}^{2^{t-1}a_n}\overline{F}(y)dy\Big)\nonumber\\
&\leq&4bt^{-1}+O(a_n^{2\alpha-1})<\infty.
\end{eqnarray*}
Recall that, for two positive functions $f$ and $g$, notation $f(x) \lesssim g(x)$ means that $\limsup_{x\to\infty} f(x)/g(x) \le 1$.

In the case (iii), representations (\ref{exam201-5}) and (\ref{exam501-3}) show that
\begin{eqnarray*}
T(x)/\overline{H}(x)&=&\int_{x/2}^{2^{t+1}a_n}\overline{H}(x-y)h(y)dy/\overline{H}(x)\nonumber\\
&\lesssim&bt^{-1}\int_{2^{t}a_n}^{2^{t+1}a_n}\overline{H}(2^{t+1}a_n-y)\overline F(y-2^{t}a_n)dy
\nonumber\\&=&bt^{-1}\Big(\int_{0}^{2^{t-1}a_n}+\int_{2^{t-1}a_n}^{2^{t}a_n}\Big)\overline{H}(y)\overline F(2^{t}a_n-y)dy\nonumber\\&\leq&bt^{-1}\overline F(2^{t-1}a_n)\int_{0}^{2^{t-1}a_n}\overline{H}(y)dy
+bt^{-1}\overline{H}(2^{t-1}a_n)\int_{0}^{2^{t-1}a_n}\overline F(y)dy\nonumber\\
&=&O(a_n^{2\alpha-1})<\infty.
\end{eqnarray*}

We may conclude that (\ref{309}) holds, therefore $F^{*2}\in\mathcal{OS}$.

Finally, since
$F\notin\mathcal{L}$ and $F^{*2}\in\mathcal{L}$, Proposition \ref{pron1} leads
to the conclusion that $F\notin\mathcal{OS}$.
\hfill$\Box$

\vspace{0.2cm}

{\bf Proof of (2).} Since $F\notin\mathcal{L}$, we have $F^{*\tau}\notin\mathcal{S}$,
by \cite{EmbGolVer1979}.
Under condition (\ref{thm101}), $F^{*\tau}\in\mathcal{L}$ follows from $F^{*2}\in\mathcal{L}$ and part (3) of Theorem \ref{th1}.

Under condition (\ref{thm102}) with any fixed $0<\varepsilon<1$ and $M=M(\varepsilon)\ge n$ large enough, Corollary \ref{cor001} implies that
$$(1-\varepsilon)\overline{F^{*\tau}}(x-1)\le\sum_{n=1}^M p_n\overline{F^{*n}}(x-1)\le(1+\varepsilon)\sum_{n=1}^M p_n\overline{F^{*n}}(x)\le(1+\varepsilon)\overline{F^{*\tau}}(x),$$
for $x$ large enough. Since $\varepsilon >0$ is arbitrary, we get
$F^{*\tau}\in\mathcal{L}$.

Further, we prove that $F^{*\tau}\in\mathcal{OS}$ under condition (\ref{thm102}). Without loss of generality, we may assume that $p_M>0$. By $F^{*2}\in\mathcal{OS}$ and Proposition 2.6 in \cite{ShiWat2005}, we have $F^{*M}\in\mathcal{OS}$. Further, by (\ref{thm102}), we have
\begin{eqnarray*}
(1-\varepsilon)\overline{F^{*\tau}}(x)&\le&\sum_{i=1}^Mp_i\overline{F^{*i}}(x)=O(\overline{F^{*M}}(x)).
\end{eqnarray*}
On the other hand, relation $\overline{F^{*M}}(x)=O(\overline{F^{*\tau}}(x))$ is clear.
Therefore, $F^{*\tau}\in\mathcal{OS}$ follows from $F^{*M}\in\mathcal{OS}$.

Next, we prove (\ref{thm1020}). Recall that all distributions from the class
$\mathcal{F}_1(0)$ are supported by the positive half-line.
Since ${\mathbf E}X_1=\infty$ and
${\mathbf E}\tau<\infty$, Theorem 1 of Denisov et al. \cite{DenFosKor2008} implies the first equality in \eqref{thm1020} (see also \cite{Rud1973}, for a particular case of power
tails). Then the first inequality in \eqref{thm1020}
follows, say, by \cite{FosKor2007}. Further, since $\tau \ge 2 [\tau /2]$ a.s.
(here $[x]$ is the integer part of $x$),
the second inequality is straightforward:
\begin{eqnarray*}
\liminf\overline{F^{*\tau}}(x)/\overline{F^{*2}}(x)
&\ge &
\liminf\overline{F^{*2[\tau /2]}}(x)/\overline{F^{*2}}(x)\\
&= &
\liminf\sum_{m=1}^\infty(p_{2m}+p_{2m+1})\overline{F^{*2m}}(x)/\overline{F^{*2}}(x)\\
&\ge &
\sum_{m=1}^\infty(p_{2m}+p_{2m+1})\liminf \overline{F^{*2m}}(x)/\overline{F^{*2}}(x)
=\sum_{m=1}^\infty(p_{2m}+p_{2m+1})m,
\end{eqnarray*}
where the last equality follows again by \cite{DenFosKor2008}.

Finally, we prove (\ref{thm1021}). Since $F\notin\cal{OS}$ and $F$ is supported by the
positive half-line, the last equality in
\eqref{thm1021} follows. By $F^{*2}\in\mathcal{L}\cap\mathcal{OS}$ and the
corresponding Kesten$^{,}$s type inequality, see Lemma 5 in \cite{YW2014}, for any $\varepsilon>0$ there is a constant $K=K(\varepsilon)>0$ such that, for all $n\ge1$ and $x\ge0$,
$$\overline{F^{*2m}}(x)/\overline{F^{*2}}(x)\le
 K(C^*(F^{*2})-1+\varepsilon )^m.$$
Further, by Lemma 4 or Remark 2 in \cite{YW2014}, for all $m\ge1$,
$$\limsup\overline{F^{*2m}}(x)/\overline{F^{*2}}(x)\le m(C^*(F^{*2})-1)^{m-1}.$$
Thus, by condition (\ref{rem104}) with $n=2$ and the dominated convergence theorem, we obtain the first inequality in (\ref{thm1021}):
\begin{eqnarray*}
\limsup\overline{F^{*\tau}}(x)/\overline{F^{*2}}(x)
&\le& \limsup\overline{F^{*2[(\tau +1)/2]}}(x)/\overline{F^{*2}}(x)\\
&=&\limsup\sum_{m=1}^\infty(p_{2m-1}+p_{2m})\overline{F^{*2m}}(x)/\overline{F^{*2}}(x)\nonumber\\
&\le&\sum_{m=1}^\infty m(p_{2m-1}+p_{2m})(C^*(F^{*2})-1)^{m-1}<\infty.
\end{eqnarray*}

\hfill$\Box$

{\bf Proof of (3).} Let $H$ be an infinitely divisible distribution on the positive half-line.
The Laplace transform of $H$ is given by
$$
\int_0^\infty \exp\{-\lambda
y\}H(dy)=\exp\{-a\lambda-\int_0^\infty(1-e^{\lambda y})\upsilon(dy)\}
$$
where $a\ge0$ is a constant and the L$\acute{e}$vy
measure $\upsilon$ is a Borel measure supported by $(0,\infty)$ with the
properties $\mu=\upsilon((1,\infty))<\infty$ and $\int_0^1y\upsilon(dy)
<\infty$ -- see, for example, Feller \cite{Fel1971}, page 450. Let
$F(x)=\mu^{-1}\upsilon(x)=\mu^{-1}\upsilon((0,x])$ for $x>0$.

It is well-known that the distribution $H$ admits the
representation $H=H^{(1)}*H^{(2)}$, where
$\overline {H^{(1)}}(x)=O(e^{-\beta x})$ for
some $\beta>0$ and
$$
H^{(2)}(x)=e^{-\mu}\sum_{n=0}^\infty\frac{\mu^n}{n!}F^{*n}(x).
$$

Let a random variable $\tau$ have a Poisson distribution,
$p_n=e^{-\mu}\frac{\mu^n}{n!}$ for $n=0,1,\cdots$.
Take a distribution
$F\in\mathcal{F}_1(0)$.
Since a Poisson distribution has unbounded support
and is light-tailed, condition (\ref{rem104}) is fulfilled and
$H^2\in\mathcal{L}\cap\mathcal{OS}$, by part (2) of Theorem \ref{th1}.
Since $H^{(1)}$ is light-tailed, we have $\overline{H}(x)\sim \overline{H^{(2)}}(x)$,
by Property \ref{property202}. Then, clearly, $H\in\mathcal{L}\cap\mathcal{OS}$.
Since distribution $G$ is Poisson, condition \eqref{thm1021} holds.
Finally,  since $F\notin\mathcal{S}$,
Theorem 1 of Embrechts et al. \cite{EmbGolVer1979} leads to $H\notin\mathcal{S}$.
\hfill$\Box$

\vspace{0.2cm}
\noindent{\bf Proof of Proposition \ref{pron1}.} By Theorem 3.1 (b) of Embrechts and Goldie \cite{EmbGolVer1979}, we
need to prove the implication $\Leftarrow$ only.
By $F^{*2}\in\mathcal{L}$ and Corollary \ref{cor001} (2),
we know that $G_2=:pF+qF^{*2}\in\mathcal{L}$ for any $p+q=1$ and $0<q<1$. Further,
since $F\in\mathcal{OS}$, we have $\overline{G_2}(x)=O(\overline{F}(x))$.
Therefore, $F\in\mathcal{L}$ follows from Lemma 2.4 of Yu et al. \cite{YWY2010}.\hfill$\Box$

\vspace{0.2cm}

\noindent{\bf Proof of Lemma \ref{gengen}.}
By $F^{*n}\in\mathcal{L}\cap\mathcal{OS}$ and Lemma 5 of Yu and Wang \cite{YW2014}, for any $0<\varepsilon_0<1$, there exists a constant $K=K(\varepsilon_0)>0$ such that, for all $x>0$ ang $m\ge1$,
$$
\overline{F^{*mn}}(x)\le K(C^*(F^{*n})-1+\varepsilon_0)^m\overline{F^{*n}}(x).
$$
Then, by (\ref{rem104}), for any $0<\varepsilon<1$, there exists an integer $M_0=M_0(\varepsilon)>1$ large enough such that
\begin{eqnarray*}
\sum_{k=(M_0-1)n}^\infty p_k\overline{F^{*k}}(x)\le\sum_{m=M}^\infty \Big(\sum_{k=(m-1)n+1}^{mn}p_k\Big)\overline{F^{*mn}}(x)\le\varepsilon\overline{F^{*\tau}}(x).
\end{eqnarray*}
Take $M=(M_0-1)n$, then (\ref{thm102}) holds.
\vspace{0.2cm}
\hfill$\Box$

\vspace{0.2cm}

{\bf Acknowledgement.}
The authors thank Dima Korshunov for bringing
our attention to the Embrechts-Goldie's conjecture and to the related problem for
infinitely divisible distributions, and for influential comments.
We thank Paul Embrechts for a number of useful remarks and for the additional bibliographical data. We also thank Takaaki Shimura for his translation and sending to us the  English version of Report \cite{ShiWat2005b}, and Stan Zachary for his comments on the Introductory part of the paper.
Finally, the authors are most grateful to two referees for
their valuable comments and suggestions.

\section{Appendix}
\setcounter{thm}{0}\setcounter{Corol}{0}\setcounter{lemma}{0}\setcounter{pron}{0}
\setcounter{remark}{0}\setcounter{exam}{0}\setcounter{equation}{0}

\subsection{On condition (\ref{101})}

The following two examples show the feasibility of condition (\ref{101}).

\begin{exam}\label{exam501}
Take a distribution $G_1$ given by
$$\overline{G_1}(x)=\textbf{\emph{I}}(x<0)+e^{-\sqrt{x}}\textbf{\emph{I}}(x\ge0).$$
Xu et al. \cite{XSWC2014} in their Example 2.1 introduce
a distribution $F_1$ on the positive half-line
such that, for $x\in(-\infty,\infty)$,
\begin{eqnarray*}
\overline{F_1}(x)&=&\overline{G}_1(x)\textbf{\emph{I}}(x<x_1)+
\sum_{n=1}^{\infty}\Big(\overline{G}_1(x_n)\textbf{\emph{I}}(x_n\le x<y_n)+\overline{G}_1(x)\textbf{\emph{I}}(y_n\le x<x_{n+1})\Big),
\end{eqnarray*}
where $\{x_n,n\ge1\}$ and $\{y_n, n\ge 1\}$ are two sequences of positive constants satisfying $x_n<y_n<x_{n+1}$ and $\overline{G}_1(x_n)=2\overline{G}_1(y_n),\ n\ge1$. One
can easily verify that
$F_1\in\mathcal{OL}\setminus\mathcal{L}$ and $\overline{F_1}(x)\asymp\overline{G_1}(x)$,
that is $0< \liminf \overline{G}_1(x)/\overline{F}_1(x)
\le \limsup \overline{G}_1(x)/\overline{F}_1(x)<\infty.$
Further, take a distribution $F_2$ such that
$$
\overline{F_2}(x)=\textbf{\emph{I}}(x<0)+\overline{G_1}(x)\textbf{\emph{I}}(x\ge0)/\log(x+2).
$$
Clearly, $F_2\in\mathcal{S}\subset\mathcal{L}$, $\overline{F_2}(x)=o(\overline{F_1}(x))$ and condition (\ref{101}) holds. Then Remark \ref{eq-eq} or, equivalently, part (1b)
of Theorem \ref{th1} imply that
 $F_1*F_2\in\mathcal{L}$.\hfill$\Box$
\end{exam}

\begin{exam}\label{exam503}
Assume $\overline{F_2}(x)=x^{-\alpha}$ for $x\ge 1$, where
$\alpha >0$. Let $1> \varepsilon_n\downarrow 0$ be any decreasing sequence. Given
two sequences $\{a_n,n\ge1\}$ and $\{b_n,n\ge1\}$ such that
$$
1=a_1<b_1<\ldots< a_n<b_n<a_{n+1}<b_{n+1}<\ldots,
$$
we let
$$
\overline{F_1}(x)=\textbf{\emph{I}}(x<a_1)+\sum\limits_{n=1}^{\infty}c_n
\textbf{\emph{I}}(x\in[a_n,b_n])
+\sum\limits_{n=1}^{\infty}d_nx^{-2\alpha}\textbf{\emph{I}}(x\in(b_n,
a_{n+1})).
$$
Here $c_1=1$, $d_n=c_nb_n^{2\alpha}$ and $c_{n+1}=d_na_{n+1}^{-2\alpha}\varepsilon_{n}.$
Then we may determine sequences $\{a_n,n\ge1\}$ and $\{b_n,n\ge1\}$ recursively in such a way that
\begin{equation}\label{anbn}
\frac{\overline{F_1}(b_n)}{\overline{F_2}(b_n)}=c_nb_n^{\alpha}=2^n\rightarrow\infty \quad \mbox{and}
\quad
\frac{\overline{F_1}(a_n-0)}{\overline{F_2}(a_n)}=d_{n-1}a_n^{-\alpha}=2^{-n+1}\rightarrow0, \quad\ \ as\ \ n\to\infty.
\end{equation}

Informally, we proceed as follows. Let $a_1=c_1=1$ and choose $b_1$ such that
$\overline{F_2}(b_1)=1/2$, then $d_1=b_1^{2\alpha}$. Then choose $a_2$ such that
$d_1a_2^{-2\alpha}=2^{-1}a_2^{-\alpha}$ and then
$c_{2}=d_1a_2^{-2\alpha}\varepsilon_1$. By the induction argument, given
$a_n$ and $c_n$, we keep $F_1(x)$ constant in the interval $[a_n,b_n]$. Since
$\overline{G}$ decreases to 0 continuously, we may choose $b_n$ so large
that the first equation in \eqref{anbn} holds. Then, by the symmetric argument,
we may choose $a_{n+1}$ so large that the second equation in \eqref{anbn} holds,
with $a_{n+1}$ in place of $a_n$.

One can see that $F_1\notin \mathcal{OL}$. However, condition (\ref{101}) is satisfied, thus $F=F_1*F_2\in\mathcal{L}$, by
Remark \ref{eq-eq} or
part (1b) of Theorem \ref{th1}.
\hfill$\Box$
\end{exam}

%
%

\subsection{Sketch of the Proof of Proposition \ref{classF2}}

The proof mostly follows the lines of the proof of Theorem \ref{th2}, so we provide
its sketch only, and also a complete proof of the last new statement.

We first analyse the distribution $F$ of the random variable $\xi$ and its density $f$.
Clearly, $\eta^{1/t}\leq\xi\leq(2\eta)^{1/t}$. Since
$$
\textbf{E}\eta^{s/t}=C\sum_{n=0}^{\infty}a_n^{s/t-\alpha}<\infty
$$
if and only if $s < t\alpha$, the same holds for $\xi$, and the distribution
$F$ is heavy-tailed with infinite mean.
Next, for all $x$, we get
\begin{eqnarray*}
&&\textbf{P}(a_n^{1/t}<\xi\le
x)\textbf{I}(x\in[a_n^{1/t},a_{n+1}^{1/t}))\\
&=&C
a_n^{-\alpha}\lo(a_n^{-1}x^{t}-1\ro)\textbf{I}(x\in[a_n^{1/t},(2a_n)^{1/t}))
+C a_n^{-\alpha}\textbf{I}(x\in[(2a_n)^{1/t},a_{n+1}^{1/t})),
\end{eqnarray*}
then
\begin{eqnarray*}
f(x)=Ct \sum_{n=0}^{\infty}x^{t-1}a_n^{-\alpha-1}\textbf{I}(x\in
[a_n^{1/t},(2a_n)^{1/t}))
\end{eqnarray*}
and
\begin{eqnarray*}
\overline{F}(x)&=&\textbf{I}(x< a_0)
+C\sum\limits_{n=0}^{\infty}\Big(\big(\sum\limits_{i=n}^{\infty}
a_i^{-\alpha}-a_n^{-\alpha}(a_n^{-1} x^t-1)\big)\textbf{I}(x\in
[a_n^{t^{-1}},(2a_n)^{t^{-1}}))\nonumber\\
&&+\sum\limits_{i=n+1}^{\infty}a_i^{-\alpha}\textbf{I}(x\in[(2a_n)^{1/t},a_{n+1}^{1/t}))\Big).
\end{eqnarray*}

Then we follow the lines of the proof of Theorem \ref{th2} to show that
$F\notin\mathcal{OL}$ and that $F^{*2}\in\mathcal{L}$, by considering again
the three cases.

Then we come to the proof of the two last statements:
 $F^{*2}\in\mathcal{OS}$ if $\alpha\in[1/2,1/t)$, and
$F^{*2}\notin\mathcal{OS}$ if $\alpha\in(1-1/t,1/2)$.
The proof in the case
$\alpha\in[1/2,1/t)$ is again analogous to the corresponding part of the proof
of Theorem \ref{th2},
so we turn to the proof of the latter result.

Let again $H=F^{*2}$, $h$ be the density of $H$, and $
T(x)=\int_{x/2}^x \overline{H}(x-y)h(y)dy$. For $\alpha\in(1-1/t,1/2)$, we have
\begin{eqnarray*}
&&T(2(2a_n)^{1/t})/\overline{H}(2(2a_n)^{1/t})\nonumber\\
&\geq&\int_{(2^{1/t}+1)a_n^{1/t}}^{2(2a_n)^{1/t}}2\overline{H}(2(2a_n)^{1/t}-y)
\int_{y/2}^{(2a_n)^{1/t}}f(y-z)f(z)dzdy/\overline{H}(2(2a_n)^{1/t})\nonumber\\
&\gtrsim&ta_n^{1-1/t}\int_{(2^{1/t}+1)a_n^{1/t}}^{2(2a_n)^{1/t}}\overline{H}(2(2a_n)^{1/t}-y)(\overline F(y-(2a_n)^{1/t})-\overline F(y/2))dy\nonumber\\
&\geq&Cta_n^{-\alpha-1/t}\int_{(2^{1/t}+1)a_n^{1/t}}^{2(2a_n)^{1/t}}\overline{H}(2(2a_n)^{1/t}-y)
((2a_n)^{1/t}-y/2)(y/2)^{t-1}dy\nonumber\\
&\geq&C2^{-1}ta_n^{1-\alpha-2/t}\int_{0}^{(2^{1/t}-1)a_n^{1/t}}\overline{H}(y)ydy\to\infty,~~~n\to\infty.
\end{eqnarray*}
Here notation $f(x)\gtrsim g(x)$ 
is equivalent to $g(x)\lesssim f(x)$ and means that $\liminf_{x\to\infty} f(x)/g(x)\ge 1$.
Thus, $F^{*2}\notin\mathcal{OS}$.

\end{document}